\documentclass{article}%
\usepackage{amsmath}%
\usepackage{amsfonts}%
\usepackage{amssymb}%
\usepackage{graphicx}

\begin{document}

\title{On Infinite Product Identities Generating Solutions for Series}
\author{Henrik Stenlund\thanks{The author is obliged to Visilab Signal Technologies for supporting this work.}
\\Visilab Signal Technologies Oy, Finland}
\date{30th October, 2016}
\maketitle
\begin{abstract}
In this paper we present a new identity and some of its variants which can be used for finding solutions while solving fractional infinite and finite series. We introduce another simple identity which is capable of generating solutions for some finite series. We demonstrate a method for generation of variants of the identities based on the findings. The identities are useful for solving various infinite products. \footnote{Visilab Report \#2016-10}\footnote{Typos removed from equations 20, 21}

\subsection{Keywords}
Summation of finite series, infinite series, identities, infinite products
\subsection{Mathematical Classification}
MSC: 20F14, 40A25, 65B10 40A20
\end{abstract}
\tableofcontents
\section{Introduction}
\subsection{General}
The motivation for this paper has been the need to treat fractional infinite and finite series of the form below,
\begin{equation}
A_N=\sum_{k=1}^{N}{\frac{ak+b}{(ck+d)(ek+f)}} \label{eqn2}
\end{equation}
or
\begin{equation}
A_N=\sum_{k=1}^{N}{\frac{(ak+b)(gk+h)}{(ck+d)(ek+f)(lk+m)}} \label{eqn4}
\end{equation}
or even higher order than these and with other functional content. The $a$ to $m\in{C}$ are constant coefficients. In this paper we show a few cases like this which can be solved, some with $N$ infinite, some finite. 

Not many methods exist for tackling series like these. Differentiation or integration with respect to a parameter may in some cases lead to a solution. Other similar tricks may be utilized. Some particular series have been tabulated, for instance in \cite{Gradshteyn2007} and \cite{Jeffrey2008}.
\subsection{Identities}
Simple identities are well-known but are, after all, far too trivial to be useful for forming bases for summation rules. They tend to boil down when written out. A parameter should be available, a bare algebraic structure is usually not sufficient. Therefore, we need a more complex identity with the aforementioned features. It would be nice to find some method for finding new identities of some form.

First, we present an infinite product identity and create variants of it. Then we generate the main summation formulas from each. Appendix A shows in detail how the basic identity is arrived at. We also write down a very simple identity which might be at the upper boundary of triviality in order to get some summation formula. Then we show the method for deriving new infinite product identities and summation rules thereof.
\section{The Identity and Its Variants}
\subsection{The Identity}
The identity in the following is derived in Appendix A and we take the result here as given. 
\begin{equation}
s=\prod_{k=2}^{\infty}{\frac{k(s(k-1)+2)}{(k+1)(s(k-2)+2)}} \label{eqn100}
\end{equation}
with $s\in{C}$, both sides being analytic. This identity is a bit strange since it is a nontrivial mapping of the variable $s$ to itself. 
\subsection{A Parametrized Identity}
In equation (\ref{eqn100}) we can spot the number $2$. The scaling properties of our variable $s$ point to the fact that we can replace it with a new parameter $\alpha$. 
\begin{equation}
s=\frac{\alpha}{2}\prod_{k=2}^{\infty}{\frac{k(s(k-1)+\alpha)}{(k+1)(s(k-2)+\alpha)}} \label{eqn150}
\end{equation}
with $s,\alpha\in{C}$ but $\alpha\neq{0}$. Again, the singularities are apparent only. This identity has a wider field of application compared to the first one (\ref{eqn100}).
\subsection{The Identity Raised to a Power}
The identity (\ref{eqn150}) can be raised to a power $r$, both sides. On the other hand, we can substitute $s$ raised to a power to the identity (\ref{eqn150}). We have two equations with the same left sides and we can then solve out the $\alpha$ with a result
\begin{equation}
(\frac{2}{\alpha})^{r-1}=\prod_{k=2}^{\infty}{\frac{k^{r-1}(s(k-1)+\alpha)^{r}(s^r(k-2)+\alpha)}{(k+1)^{r-1}(s^r(k-1)+\alpha)(s(k-2)+\alpha)^r}} \label{eqn160}
\end{equation}
with $s,\alpha\in{C}, r\in{R}$, $s,\alpha\neq{0}$
\subsection{A Parametrized Finite Identity}
It is useful to have the finite case of the parametrized identity (\ref{eqn150}). We can calculate the identity with a finite number of terms getting a new modified identity
\begin{equation}
\frac{(N-1)s+\alpha}{N+1}=\frac{\alpha}{2}\prod_{k=2}^{N}{\frac{k(s(k-1)+\alpha)}{(k+1)(s(k-2)+\alpha)}} \label{eqn165}
\end{equation}
with $s,\alpha\in{C}$, $\alpha\neq{0}$. $N$ is now finite. 
\subsection{A Trivial Finite Identity}
We are able to develop a very simple identity to make it suitable for further development as a summation formula. The following equation is very easy to prove.
\begin{equation}
\frac{s}{s+N}=\prod_{j=1}^{N}{(1-\frac{1}{s+j})} \label{eqn170}
\end{equation}
with $s\in{C}$ and $N=1,2,3$ finite. Let's set
\begin{equation}
\frac{s}{s+N}=z \label{eqn180}
\end{equation}
and substitute for $s$ to end up with a new identity
\begin{equation}
z=\prod_{j=1}^{N}{\frac{(j(1-z)+z(N+1)-1)}{(j(1-z)+Nz)}} \label{eqn190}
\end{equation}
with $z\in{C}$ and $N$ finite.
\section{Summation Formulas}
We convert the identities obtained in the preceding paragraphs to infinite or finite sums applicable as summation formulas. We use a similar technique in each case, i.e. taking a logarithm and then differentiating the expressions with respect to some parameter. Direct differentiation or integration or other algebraic operations may sometimes deliver useful summation relations.
\subsection{The Basic Summation Formula}
We start by taking the logarithm of the identity (\ref{eqn100}) above
\begin{equation}
ln(s)=\sum_{k=2}^{\infty}{ln[\frac{k(s(k-1)+2)}{(k+1)(s(k-2)+2)}]} \label{eqn310}
\end{equation}
This can be differentiated with respect to $s$ to yield
\begin{equation}
\frac{1}{s}=2\sum_{k=2}^{\infty}{\frac{1}{(s(k-1)+2)(s(k-2)+2)}} \label{eqn320}
\end{equation}
still with $s\in{C}$, $s\neq{0}$.
\subsection{The Parametrized Summation}
Another summation formula is obtained from identity (\ref{eqn150}) by taking the logarithm and differentiating it with respect to $\alpha$ to arrive at
\begin{equation}
\frac{1}{\alpha}=\sum_{k=2}^{\infty}{\frac{s}{(s(k-1)+\alpha)(s(k-2)+\alpha)}} \label{eqn500}
\end{equation}
with $s,\alpha\in{C}$, $s,\alpha\neq{0}$. We do have here two parameters available for fitting to summation problems.
\subsection{Summation of the Power Expression}
If we take logarithms of both sides of (\ref{eqn160}) and differentiate with respect to $\alpha$ we obtain
\begin{equation}
\frac{1-r}{\alpha}=\sum_{k=2}^{\infty}{\frac{-rs}{(s(k-1)+\alpha)(s(k-2)+\alpha)}+\frac{s^r}{(s^r(k-1)+\alpha)(s^s(k-2)+\alpha)}} \label{eqn600}
\end{equation}
with three free parameters $s,\alpha\in{C}, r\in{R}$, $s,\alpha\neq{0}$. The poles are apparent again.
\subsection{The Parametrized Finite Summation}
Again, we subject the expression (\ref{eqn165}) to logarithm and differentiation with respect to $\alpha$
\begin{equation}
\frac{1}{\alpha}-\frac{1}{(N-1)s+\alpha}=s\sum_{k=2}^{N}{\frac{1}{(s(k-1)+\alpha)(s(k-2)+\alpha)}} \label{eqn700}
\end{equation}
with $s,\alpha\in{C}$, $s,\alpha\neq{0}$. We have two parameters available for fitting to summation problems, in addition to the number of terms. This is a known result, see \cite{Gradshteyn2007}. We could do this same procedure for all infinite summations presented here getting finite summation formulas.
\subsection{The Trivial Summation}
We process the last identity (\ref{eqn190}) by taking the logarithm and differentiating it with respect to $z$
\begin{equation}
\frac{1}{z}=N\sum_{j=1}^{N}{\frac{1}{[j(1-z)+z(N+1)-1][j(1-z)+Nz]}} \label{eqn800}
\end{equation}
with $z\in{C}, z\neq{0}$ and $N=1,2,3...$ finite. We can also attack the trivial identity (\ref{eqn170}) directly by differentiating it with respect to $s$ getting
\begin{equation}
\frac{N}{(s+N)s}=\sum_{k=1}^{N}{\frac{1}{(s+k)(s+k-1)}} \label{eqn850}
\end{equation}
This result happens to be already known, related to equation (\ref{eqn700}).
\section{Synthesizing New Identities}
We soon realize that  we can produce new identities as we have the idea in front of us. We can add any number of terms to the product with balancing terms in the denominator. We can add terms to the product of various powers and combinations. We can adjust the index offsets causing unexpected effects on the final result in the form of functional divisors.
\subsection{Creating a Power Identity}
As an example let's substitute $s^r$ to our identity and raise some of the terms to the same power. Thus we have
\begin{equation}
s^{r}=\frac{\alpha}{2^r}\prod_{k=2}^{\infty}{\frac{k^{r}(s^r(k-1)^r+\alpha)}{(k+1)^{r}(s^r(k-2)^r+\alpha)}} \label{eqn1200}
\end{equation}
with $s,\alpha\in{C}, r\in{R}$, $s,\alpha\neq{0}$
The corresponding summation is derived as usual to become
\begin{equation}
\frac{1}{s}=\alpha{s^{r-1}}\sum_{k=2}^{\infty}{\frac{(k-1)^r-(k-2)^r}{(s^r(k-1)^r+\alpha)(s^r(k-2)^r+\alpha)}} \label{eqn1250}
\end{equation}
It is obvious that we can proceed a long way in this direction and generate an infinity of summation formulas of fractional type. 
\subsection{A Method}
We have now an idea of how to create simple infinite products by using equation (\ref{eqn150}). We can generalize it to the next level by using a function $f(x)$
\begin{equation}
\beta\prod_{k=2}^{\infty}{\frac{f(k)(f(s(k-1))+\alpha)}{f(k+1)(f(s(k-2))+\alpha)}} \label{eqn2000}
\end{equation}
with $s,\alpha,\beta\in{C}$. $\beta$ is a scaling parameter to be used when we have evaluated the product. The product will come out in the form of
\begin{equation}
\beta\frac{f(2)(f(s(N-1))+\alpha)}{(f(0)+\alpha)f(N+1)} \label{eqn2020}
\end{equation}
or if we let $N\rightarrow\infty$
\begin{equation}
\beta\frac{f(2)f(sN)}{(f(0)+\alpha)f(N)} \label{eqn2030}
\end{equation}
In order to get a finite nonzero value of the product we have at least the following options for the $f(x)$,
\begin{equation}
\lim_{N\rightarrow\infty}\frac{f(sN)}{f(N)}=C_0s^r \label{eqn2032}
\end{equation}
or the $f(x)$ saturates to a limiting value while $x\rightarrow\infty$. $f(x)$ must not be oscillating with a sign change nor grow in another way than indicated in equation (\ref{eqn2032}). The sum is constructed as in preceding paragraphs as soon as the new product is evaluated. Convergence of the new product should be checked if not obvious.
\subsubsection{A Fractional Function $f(x)$}
As a first example we select 
\begin{equation}
f(x)=\frac{x^2}{1+x} \label{eqn2080}
\end{equation}
This function will asymptotically approach $x$ while $x\rightarrow\infty$. Then we get the result
\begin{equation}
s=\frac{3\alpha}{4}\prod_{k=2}^{\infty}{\frac{(\frac{k^2}{1+k})[\frac{s^2(k-1)^2}{(1+s(k-1))}+\alpha]}{(\frac{(1+k)^2}{2+k})[\frac{s^2(k-2)^2}{(1+s(k-2))}+\alpha]}} \label{eqn2090}
\end{equation}
with $s,\alpha\in{C}, \alpha\neq{0}$ and we have already solved the product with the value on the left.
We can furnish the sum as follows
\begin{equation}
\frac{1}{s}=s\sum_{k=2}^{\infty}{[\frac{2(k-1)^2+(k-1)^3}{s^2(k-1)^2+s^3(k-1)^3+\alpha(1+s(k-1))^2}} \nonumber
\end{equation}
\begin{equation}
-{\frac{2(k-2)^2+(k-2)^3}{s^2(k-2)^2+s^3(k-2)^3+\alpha(1+s(k-2))^2}}] \label{eqn2220}
\end{equation}
\subsubsection{A $tanh(x)$ as a Function}
As a second example we take
\begin{equation}
f(x)={tanh(x)} \label{eqn2230}
\end{equation}
This function will asymptotically approach a constant while $x\rightarrow\infty$. After substituting this to (\ref{eqn2000}) we have the product
\begin{equation}
1=\frac{\alpha}{(1+\alpha)tanh(2)}\prod_{k=2}^{\infty}{\frac{{tanh(k)}[tanh(s(k-1))+\alpha]}{{tanh(k+1)}[tanh(s(k-2))+\alpha]}} \label{eqn2240}
\end{equation}
with $s,\alpha\in{C}, s,\alpha\neq{0}$. Note that the value of this product is independent on $s$.
The summation is equal to
\begin{equation}
\frac{1}{\alpha(1+\alpha)}=\sum_{k=2}^{\infty}{\frac{1}{tanh(s(k-2))+\alpha}-\frac{1}{tanh(s(k-1))+\alpha}} \label{eqn2250}
\end{equation}
\subsubsection{A Function of form $x\cdot{tanh(x)}$}
As a third example we select
\begin{equation}
f(x)=x\cdot{tanh(x)} \label{eqn2280}
\end{equation}
This function will asymptotically approach $x$ while $x\rightarrow\infty$. We substitute this to (\ref{eqn2000}) and get
\begin{equation}
s=\frac{\alpha}{2tanh(2)}\prod_{k=2}^{\infty}{\frac{k\cdot{tanh(k)}[s(k-1)tanh(s(k-1))+\alpha]}{(k+1)\cdot{tanh(k+1)}[s(k-2)tanh(s(k-2))+\alpha]}} \label{eqn2290}
\end{equation}
with $s,\alpha\in{C}, s,\alpha\neq{0}$.
The summation will finalize as
\begin{equation}
\frac{1}{s}=\sum_{k=2}^{\infty}{[\frac{(k-1)tanh(s(k-1))+\frac{s(k-1)^2}{cosh^2(s(k-1))}}{s(k-1)tanh(s(k-1))+\alpha}} \nonumber
\end{equation}
\begin{equation}
-\frac{(k-2)tanh(s(k-2))+\frac{s(k-2)^2}{cosh^2(s(k-2))}}{s(k-2)tanh(s(k-2))+\alpha}] \label{eqn2420}
\end{equation}
\subsubsection{An $x\cdot{arctan(x)}$ Function}
As a fourth example we pick up
\begin{equation}
f(x)=x\cdot{arctan(x)} \label{eqn2480}
\end{equation}
taking the principal value only. This function will asymptotically approach $xC_1$ while $x\rightarrow\infty$. Then we obtain for the product the result
\begin{equation}
s=\frac{\alpha}{2arctan(2)}\prod_{k=2}^{\infty}{\frac{k\cdot{arctan(k)}[s(k-1)arctan(s(k-1))+\alpha]}{(k+1)\cdot{arctan(k+1)}[s(k-2)arctan(s(k-2))+\alpha]}}\label{eqn2490}
\end{equation}
with $s,\alpha\in{C}, s,\alpha\neq{0}$.
The resulting sum is the following
\begin{equation}
\frac{1}{s}=\sum_{k=2}^{\infty}{[\frac{(k-1)arctan(s(k-1))+s^2(k-1)^3arctan(s(k-1))+s(k-1)^2}{(1+s^2(k-1)^2)(s(k-1)arctan(s(k-1))+\alpha)}} \nonumber
\end{equation}
\begin{equation}
-\frac{(k-2)arctan(s(k-2))+s^2(k-2)^3arctan(s(k-2))+s(k-2)^2}{(1+s^2(k-2)^2)(s(k-2)arctan(s(k-2))+\alpha)}] \label{eqn2620}
\end{equation}
\subsubsection{An Exponential Function $\frac{x\cdot{e^x}}{1+e^x}$}
As a last example we present a bit more complicated case of
\begin{equation}
f(x)=\frac{x\cdot{e^x}}{1+e^x} \label{eqn3480}
\end{equation}
The exponential could have a parameter 	to expand its capabilities. This function will asymptotically approach $x$ while $x\rightarrow\infty$. The product can be evaluated to
\begin{equation}
s=\frac{\alpha(1+e^2)}{2e^2}\prod_{k=2}^{\infty}{\frac{\frac{ke^k}{1+e^k}[\frac{s(k-1)e^{s(k-1)}}{1+e^{s(k-1)}}+\alpha]}{\frac{(k+1)e^{k+1}}{1+e^{k+1}}[\frac{s(k-2)e^{s(k-2)}}{1+e^{s(k-2)}}+\alpha]}}\label{eqn3490}
\end{equation}
with $s,\alpha\in{C}, \alpha\neq{0}$.
The corresponding sum will be the following
\begin{equation}
\frac{1}{s}=\sum_{k=2}^{\infty}{e^{sk}[\frac{(k-1)e^{-s}+s(k-1)^2e^{-s}-\frac{s(k-1)^2e^{s(k-2)}}{1+e^{s(k-1)}}}{(\alpha+s(k-1))e^{s(k-1)}+\alpha}} \nonumber
\end{equation}
\begin{equation}
-\frac{(k-2)e^{-2s}+s(k-2)^2e^{-2s}-\frac{s(k-2)^2e^{s(k-4)}}{1+e^{s(k-2)}}}{(\alpha+s(k-2))e^{s(k-2)}+\alpha}] \label{eqn3620}
\end{equation}
proving that we are able to derive rather complex summation relations with this method.
\section{Discussion}
We have presented new identities as starting points for generating solutions to certain categories of infinite and finite series. The identities are either infinite or finite products of fractional type, thus far being easy for evaluation. 

We showed a method on how to construct new identities as infinite products. The process is fairly straightforward. The method presented is limited in capability but powerful for some types of products, some even very complicated. One may then be able to derive a summation formula from it. Very complex summation rules are now possible. Using this method backwards, i.e. by starting from the sum is much more requiring and no result is guaranteed for identifying any infinite product belonging to it, thus possibly allowing its evaluation. We have presented several examples illustrating both the new identities and resulting summation formulas.  

\appendix
\section{Appendix. Derivation of the Identity}
By starting from summation of a finite series with the method given recently \cite{Stenlund2016} we get the following expression, formally
\begin{equation}
\sum_{j=0}^{N-1}{\frac{1}{j+a}}=\sum_{k=1}^{\infty}{(-1)^{k+1}\zeta(k+1)[(a+N-1)^k-(a-1)^{k}]} \label{eqn5000}
\end{equation}
The function is the Riemann zeta. It is quickly recognized that the expression will fail in general. However, it is valid for $N=1$
 and it will work as our starting point. We substitute $b$ for $\frac{1}{a}$ to have
\begin{equation}
b=\sum_{k=1}^{\infty}{\frac{(-1)^{k+1}\zeta(k+1)}{b^{k}}}+\sum_{k=1}^{\infty}{\frac{(-1)^{k}\zeta(k+1)(1-b)^k}{b^{k}}} \label{eqn5100}
\end{equation}
We can apply the known equation below to both of the sums above
\begin{equation}
\sum_{k=1}^{\infty}{(-1)^{k+1}s^{k+1}\zeta(k+1)}=\sum_{k=1}^{\infty}{\frac{s^2}{k(k+s)}} \label{eqn5200}
\end{equation}
yielding
\begin{equation}
\frac{1}{b}=\sum_{k=1}^{\infty}{\frac{1}{(bk+1)(bk+1-b)}} \label{eqn5300}
\end{equation}
This equation has poles at $b=\frac{-1}{N}$ but will otherwise converge for $b>0$. We integrate this expression with respect to $b$ from $1$ to $b$ and clean up the result to
\begin{equation}
ln(\frac{2b}{b+1})=\sum_{k=2}^{\infty}{ln(\frac{k(bk+1)}{(bk-b+1)(k+1)})} \label{eqn5400}
\end{equation}
We put this expression as an argument of the exponential function with the result
\begin{equation}
\frac{2b}{b+1}=\prod_{k=2}^{\infty}{\frac{k(bk+1)}{(bk-b+1)(k+1)}} \label{eqn5500}
\end{equation}
Finally, we substitute 
\begin{equation}
b=\frac{s}{2-s} \label{eqn5600}
\end{equation}
and obtain
\begin{equation}
s=\prod_{k=2}^{\infty}{\frac{k(s(k-1)+2)}{(k+1)(s(k-2)+2)}} \label{eqn6500}
\end{equation}
with $s\in{C}$. The denominator has singularities of pole type. However, when looked at more closely, we can see that the poles are apparent only. This infinite product can be analyzed in a conventional way by taking the upper limit as $N$. Then one can study its behavior as $N\rightarrow{\infty}$, leaving out infinitesimal terms.
\end{document}